\documentclass[a4paper,12pt]{article}
\usepackage{arxiv}

\usepackage[utf8]{inputenc} 
\usepackage[T1]{fontenc}    
\usepackage{hyperref}       
\usepackage{url}            
\usepackage{booktabs}       
\usepackage{amsfonts}       
\usepackage{nicefrac}       
\usepackage{microtype}      
\usepackage{lipsum}
\usepackage{graphicx}
\usepackage[
    backend=biber,
    style=numeric,
    maxnames=99
]{biblatex}
\usepackage{xcolor}
\usepackage[utf8]{inputenc}
\usepackage{titlesec}
\usepackage[british]{babel}
\usepackage{amsmath,amssymb,amsfonts}
\usepackage{amsthm}
\usepackage{graphicx}
\usepackage{multicol}
\usepackage{subcaption}
\usepackage{float}
\usepackage{geometry}
\usepackage{tikz}
\usepackage{dsfont}
\usepackage{csquotes}
\usetikzlibrary{automata, positioning}
\theoremstyle{definition}

\graphicspath{ {./images/} }

\title{Counterintuitive Problems In Discrete Probability}

\author{
 Luca Avena$^1$ \\
 luca.avena@unifi.it
   \And
 Gianmarco Bet$^1$ \\
 gianmarco.bet@unifi.it\\[0.5cm]
 $^1$Università degli Studi di Firenze\\
 Dipartimento di Matematica e Informatica \enquote{Ulisse Dini}
  \And
 Bernardo Busoni$^1$ \\
 bernardo.busoni@edu.unifi.it
}

\begin{document}
\maketitle
\begin{abstract}
This manuscript contains a collection of counterintuitive problems in discrete probability, together with detailed solutions. The dataset was constructed as part of a broader research project investigating the capabilities of the latest-generation Large Language Models (LLMs) in solving discrete probability problems, in order to assess whether LLMs tend to make systematic reasoning errors associated with known cognitive biases \cite{avena2026reliablellmscomesplaying}.
The problems collected here are specifically designed to challenge heuristic reasoning strategies that often lead to intuitively appealing but mathematically incorrect conclusions. In this context, we refer to a problem as \textit{counterintuitive} when commonly used heuristic reasoning strategies tend to produce answers different from the mathematically correct solution. The dataset combines several types of problems. Some are adapted from classical probabilistic paradoxes and cognitive-bias literature, while others originate from recreational mathematics sources or were developed by ourselves following similar principles. The primary purpose of this document is to provide a transparent and publicly accessible reference for the problems used in our experimental evaluation of language models, as well as providing detailed solutions human made. At the same time, we believe that this collection may also prove useful for future research on probabilistic reasoning, cognitive biases, and the evaluation of reasoning capabilities in artificial intelligence systems.
		


\end{abstract}

\section{Introduction}
Probability is one of those areas of mathematics in which human intuition is especially prone to error. Even when a problem can be stated in very simple terms the answer suggested by immediate intuition often differs from the mathematically correct one. This discrepancy has long attracted the attention of cognitive psychologists, from the work of Kahneman and Tversky on heuristics and biases to the many subsequent studies investigating errors in probabilistic judgement and decision-making under uncertainty. But the tension between probability and intuition is also present in classical games, paradoxes, and recreational problems, many of which continue to provide striking examples of how difficult probabilistic reasoning can be.\\
This document contains a collection of counterintuitive problems in discrete probability, together with detailed solutions. The dataset was constructed as part of a broader research project investigating the capabilities of the latest-generation Large Language Models (LLMs) in solving discrete probability problems, in order to assess whether LLMs tend to make systematic reasoning errors associated with known cognitive biases \cite{avena2026reliablellmscomesplaying}.\\
The problems collected here are specifically designed to challenge heuristic reasoning strategies that often lead to intuitively appealing but mathematically incorrect conclusions. In this context, we refer to a problem as \textit{counterintuitive} when commonly used heuristic reasoning strategies tend to produce answers different from the mathematically correct solution. The dataset combines several types of problems. Some are adapted from classical probabilistic paradoxes and cognitive-bias literature, while others originate from recreational mathematics sources or were developed by ourselves following similar principles. All the problems we have taken from well-known sources have been carefully reformulated to prevent language models from recognising them and to ensure that they provide answers based on genuine, independent reasoning. Although the problems are formulated using elementary mathematical concepts, many of them require careful probabilistic reasoning and are intentionally structured to expose superficial pattern matching or intuitive reasoning failures. Most of the exercises were designed as open-ended questions. This decision was motivated by the need to measure the capabilities of the models more accurately and to avoid any potential leading effects inherent in multiple-choice questions. On the other hand, we believe that if these problems are to be used to investigate issues relating to human reasoning, it is more appropriate to formulate the questions as multiple-choice questions.\\
The primary purpose of this document is to provide a transparent and publicly accessible reference for the problems used in our experimental evaluation of language models \cite{avena2026reliablellmscomesplaying}, as well as providing detailed solutions human made. At the same time, we believe that this collection may also prove useful for future research on probabilistic reasoning, cognitive biases, and the evaluation of reasoning capabilities in artificial intelligence systems.\\
Feedback on possible mistakes or typos and suggested additions to the list of problems are welcome at the email address of the authors.
\section{Problems}

\subsection{Coins Under Cups \cite{litt2024_probabilitypuzzle}}\label{Es1}
Alice and Bob take part in a television programme. In the proposed game, the host tosses 100 fair coins and, without showing the results to the two contestants, hides them under 100 opaque cups numbered from 1 to 100 with the result (heads or tails) facing upwards. The host then starts the game and the two begin to check the results of the coins at the same rate: Alice checks the coins under the cups in increasing order (1,2,3,...,100), whereas Bob first checks the coins under the odd-numbered cups and then those under the even-numbered cups, (1,3,...,99,2,4,...100). Suppose that the results seen by one are not visible to the other; for example, while Alice checks the coin under cup 2, she does not see the result that Bob is looking at at that moment under cup 3, and vice versa. The first person to see two heads wins. \\[0.2cm]
What is the ratio between the probability that Alice wins and the probability that Bob wins? \\[0.2cm]
\textbf{Answer:} $\sim 0,75$

\subsection{Head-Head or Head-Tails?\cite{litt2024_probabilitypuzzle}}\label{Es2}
Alice and Bob take part in a game. A fair coin is tossed 100 times by the host. Alice scores one point whenever two consecutive heads occur, whereas Bob scores one point whenever a head is immediately followed by a tail. The player who accumulates more points wins. \\[0.2cm]
What is the ratio between the probability that Alice wins and the probability that Bob wins? \\[0.2cm]
\textbf{Answer:} $\sim 0,94$

\subsection{A Sequence Race \cite{litt2024_probabilitypuzzle}}\label{Es3}
Alice and Bob take part in a game. A fair coin is tossed until one of the following two occurrences appears: CTTTC or CTCTC, where C denotes Tails and T denotes Heads. In the first case Alice wins, in the second case Bob wins.\\[0.2cm]
What is the ratio between the probability that Alice wins and the probability that Bob wins? \\[0.2cm]
\textbf{Answer:} $\sim 1,25$

\subsection{An Unfair Coin in a Haystack}\label{Es4}
The host of a prize game throws 10,000 coins onto the table in front of Alice and Bob. He then says that 9,999 are fair, but one is unfair so that whenever it is tossed it always lands heads. At this point the host takes nine coins at random, tosses each of them 10 times and obtains heads every time. Alice says that at least one of the coins must be the biased coin, while Bob says that all the coins are fair.
\\[0.2cm] Assuming everything the host says is true, what is the probability that Alice is right?\\[0.2cm]
\textbf{Answer:} $\sim 0,48$

\subsection{The Pen of Cows and Sheep \cite{litt2024_probabilitypuzzle}}\label{Es5}
A farmer raises 150 cows and 150 sheep. One day, he selects a number at random from 0 to 150 and lets out a number of cows equal to the number selected from their pen. Then he goes to the sheep pen and lets some of them out until he removes one hundred and fifty animals in total. Finally he puts all the animals that came out into a new pen. At this point, he lets one animal at random out of the new pen, and a lovely cow comes out. \\[0.2cm]
What are the chances that another cow comes out if he lets another animal at random out of the new pen?\\[0.2cm]
\textbf{Answer:} $\frac{2}{3}$

\subsection{Theft at the Greenhouse}\label{Es6}
In a greenhouse there are 51 petunias and 50 orchids. A man steals two plants at random and put them in a bag so that they could not be seen, but Bob saw that one of the two plants was a petunia.\\[0.2cm]
On her way home Bob wonders: what are the chances that the man took two petunias?
\\[0.2cm]
\textbf{Answer:} $0,5$

\subsection{Probabilistic Escape \parencite{tversky1974judgment}}\label{Es7}
Hercules is locked in a tower that has two guards, and in order to escape he must solve the following dilemma: each of the guards tosses a fair coin 10 times without being seen and states, in order, the results he obtained. The first says he obtained heads ten times, while the second says he obtained \{T,T,C,T,C,T,T,C,C,T\}. Finally they tell him that one of the two has lied and that, if he manages to guess who it was, they will set him free. \\[0.2cm]
Assuming that the liar does not want to be found out, who is more likely to have lied?\\[0.2cm]
\textbf{Answer:} Both the answer "Equivalent" and "The second guard" are accepted.

\subsection{The Bounded Envelope Game Show \parencite{gardner1982aha}}\label{Es8}
Alice has won a prize game and must now claim her winnings. The maximum prize fund available for the evening depends on various factors; suppose that its value in a given currency is drawn uniformly at random from the set $\{10\cdot 10^6, \dots,25\cdot 10^6 \}$. Subsequently, a random integer is drawn between 0 and the available maximum and the host has placed that amount of money in one envelope and half of it in the other. At this point the host goes to Alice with the two envelopes in his hand and asks her to choose one at random. Alice takes one, opens it and finds 8 million. The host then gives her the opportunity to switch envelopes.\\[0.2cm]
What is the expected value of Alice's winnings if she chose to switch envelopes?\\[0.2cm]
\textbf{Answer:} $\sim 8,06$ million

\subsection{The Fruit Market \cite{tversky1974judgment}}\label{Es9}
At the fruit market there is an apple seller. In the apple basket there are 54 red apples and 6 green apples. A blindfolded man approaches the basket and says that he has the ability to distinguish red and green apples by touch, with an accuracy of 85\% regardless of the colour. At this point the man takes an apple at random and he affirms that it is a green apple.\\[0.2cm]
What are the chances that the apple is actually green?\\[0.2cm]
\textbf{Answer:} $\frac{17}{44}$

\subsection{Holiday at the Sea}\label{Es10}
Alice and Bob want to go to the sea when sea is calm in order to take a bath. Suppose that the sea is either calm or rough, with the same probability, and that they don't know its condition until they arrive there. Therefore, when they are on holiday, they go to the sea until they have seen more calm days than rough ones, or until their vacation leave is over. For example, if they immediately see a calm day, they stop at once.\\[0.2cm]
Assuming the holiday leave of Alice and Bob is 20 days long, what is the ratio between the average number of days with calm seas and the average number of days with rough seas that they will experience?\\[0.2cm]
\textbf{Answer:} $1$

\subsection{Some Dice Games \cite{rump2001strategies}}\label{Es11}
Alice, Eva, Bob and Carlo take part in a game. The host gives each of them a fair six-sided die. The numbers on the faces of all the dice are between 0 and 6, but their distribution is not the standard one. For example, a possible die might have the numbers 2, 2, 2, 3, 5, 5 on its faces, or 1, 1, 4, 4, 6, 6. Alice and Eva play 1000 matches. In each match, each of them rolls their own die and whoever obtains the higher result wins; if the two players roll the same number, it’s a draw. Eva wins more matches than Alice. Then Eva plays 1000 matches against Bob, but Bob wins more matches than Eva. Now Bob must play 1000 matches against Alice. Alice wins more matches than Bob. Finally Carlo plays 1000 matches against each of the other three and loses the majority of the matches against each of them. That said, the four players begin to play all together, they play 1000 matches, and each match is won by whoever obtains the highest value. With respect to this last series of matches, which of the following statements is certainly true:
\begin{enumerate}
  \item[a)] One among Alice, Eva and Bob has a greater chance of winning more matches than Carlo.
  \item[b)] Eva is more likely to win more matches than Alice.
  \item[c)] Carlo may have a greater chance of winning more matches than all the others.
  \item[d)] None of the previous options is true.
\end{enumerate}
\textbf{Answer:} c)

\subsection{A Stochastic Triathlon Race \cite{10.1111/j.2517-6161.1951.tb00088.x}}\label{Es12}
In a single running race from point A to point B, Carlo beats Anna 70\% of the time. In a single swimming race from point B to point C, Carlo beats Anna 90\% of the time. Finally, in a single cycling race from point C to point D, Carlo always wins against Anna. That said, Carlo challenges Anna in a race made up of the union of the races described in the stages: AB running, BC swimming and CD cycling. Given this information, which of the following statements is certainly true?
\begin{enumerate}
  \item[a)] Since Carlo always wins the cycling race, he will win the overall race.
  \item[b)] Given the statistics of the single races, Carlo will most probably win the overall race.
  \item[c)] Despite the statistics of the single races, Anna could have a very high probability of winning the overall race.
  \item[d)] Despite the statistics of the single races, Carlo could be in difficulty, but on average he will win slightly more overall races than Anna.
\end{enumerate}
\textbf{Answer:} c)

\subsection{The New Broker \cite{gardner1988time}}\label{Es13}
In an investment fund, each year the broker who has obtained the highest returns is awarded a prize. Last year, participating in the competition there was Alberto, who had a portfolio guaranteeing a 3\% return, and Filippo, who had a portfolio that with probability 0.49 yielded 5\% and with prob. 0.51 yielded 1\%. Alberto therefore had more chances of winning the competition; in particular, he won with probability 0.51. This year Bernardo has also been hired and participates with a portfolio structured as follows: 6\% return with prob. 0.22, 4\% return with prob. 0.22 and 2\% return with prob. 0.56. Alberto, worried, compares his own portfolio with Bernardo's, but reassures himself when he realises that, between the two of them, he has more chances of obtaining a better return, with probability 0.56. Bernardo knows this, but on the other hand, comparing his own portfolio with Filippo's, he realises that he has more chances than him.\\[0.2cm]
What is the probability that the plan of the broker who wins the fund award maximises expected returns?\\[0.2cm]
\textbf{Answer:} $\sim 0,33$

\subsection{A Huge Monopoly Board}\label{Es14}
Alice and Bob challenge each other to the following game. An annulus is divided into N sections, with N arbitrarily large. Both players own six distinct sections: Alice's are contiguous with one another and the location of the first one in chosen uniformly at random, whereas Bob's are chosen from among the remaining ones, always by drawing uniformly at random without replacement. Suppose that a token is placed on a random section of the board and that from there it is moved clockwise by as many sections as the result of rolling a fair six-sided die. Each time the token is moved to a player's section, that player scores one point. After a series of moves, the game ends at the moment when the next move would take the token to, or past, the section from which it started. The player who has scored more points at the end of the game wins.\\[0.2cm]
Supposing that Bob bets X on his own victory, in order for the bet to be fair, how much should Alice bet? (In the event of a draw, each gets back what they bet)\\[0.2cm]
\textbf{Answer:} $\sim 1,03\,$X

\subsection{New Coins on the Ground \cite{bbc2010tuesdayboy}}\label{Es15}
It is known that in the cash register of a bar there are 12 coins, of which one in four is new. A man enters the bar, buys a beer and receives two coins as change. Suppose that the coins received were given at random from those present in the cash register. As he leaves the bar, the two coins fall to the ground. Another man, entering the bar, sees the coins on the ground and says aloud: "There is at least one new coin that has fallen to the ground showing Heads." What is the probability that both coins fell showing Heads?\\[0.2cm]
\textbf{Answer:} $\frac{10}{21}$

\subsection{Chaos for the Kennel \cite{henze2018absentmindedpassengers}}\label{Es16}
In a dog training centre, 46 dogs have been trained in the following way: when released one at a time, in order, into an open area containing 46 distinct kennels, each of them is able, with probability 1, to enter its own specific kennel. If instead, for some reason, they find their own kennel occupied, they will choose to enter a randomly chosen kennel that is not already occupied. An untrained dog, together with its respective kennel, is added to the group; however, when released into the open area, it enters a random kennel among those available. Suppose that the untrained dog is released first and then, one at a time in sequence, the other 46 trained dogs. What is the probability that the last dog released enters its own kennel? \\[0.2cm]
\textbf{Answer:} $\frac{1}{2}$

\subsection{Buses on Strike}\label{Es17}
In a neighbourhood on the outskirts of Florence, a bus passes on average every 10 minutes throughout the day. Let us assume that the inter-arrival times between buses are independent and identically distributed as exponential random variables. Suppose that, a priori, the probability that on a given day the public transport service is on strike is equal to 0.01 and that each day the service is on strike independently of the other days. I am at a bus stop in the neighbourhood described; it is 10 in the morning and I have been waiting for an hour, but no bus has yet passed. \\[0.2cm]
How many minutes do I expect to have to wait still before one passes?\\[0.2cm]
\textbf{Answer:} $\sim 696$

\subsection{Everyone’s Wishes}\label{Es18}
Years ago Bob and Alice were a childless couple, and each of them had very particular wishes: Bob wanted to have at least four children, Alice wanted to have at least two male children. Moreover Alice's parents, very much wanted to have a grandchild of each sex. Over these years Bob and Alice had children until they had satisfied everyone's wishes, including Alice's parents.\\[0.2cm]
Assuming that the probabilities of a child's sex are equal and that the sexes of successive births are independent, what is the probability that the couple's second child was male?\\[0.2cm]
\textbf{Answer:} $\frac{1}{2}$

\subsection{The Dogs of Alice \cite{tversky1974judgment}}\label{Es19}
Alice is a cheerful and kind girl. She has always loved dogs and looking after them. At the dog training centre there are four dogs belonging to her, which have caused no problems in management because of how well-behaved they have been: Ulisse, a white and always cheerful Labrador, Minosse, a very large and well-groomed wolfdog, Teseo, another beautiful wolfdog who is always playful, and finally Hercules. From the data available, is it more likely that Hercules is a bad dog or that he is a bad but well-groomed dog?
\\[0.2cm]
\textbf{Answer:} That Hercules is a bad dog.

\subsection{Monty Hall, But with a Golden Sheep \cite{selvin1975problem}}\label{Es20}
We are taking part in a prize game. In front of us there are three doors; behind two of them there is a sheep, behind the last one there is a car, and we, the contestants in the game, do not know their arrangement. However, we know that one of the sheep belongs to a billionaire who lost it, was very fond of it, and would pay a great deal of money to get it back. Therefore our aim is to win the billionaire's sheep. We shall call this sheep the “golden sheep”. The host does not know about the presence of this particular sheep and therefore cannot distinguish the golden sheep from the normal sheep either. We, however, are able to distinguish the golden sheep from the other sheep. The host knows how the sheep and the car are arranged behind the doors; we do not. The game begins and we choose a door. The host opens at random one of the doors we have not chosen among those that hide a sheep behind them. We realise that the one shown is not the golden sheep. At this point the host asks us whether we want to switch the door initially chosen and win the prize behind the other still-closed door that we did not choose, or win the prize behind the first door we chose. \\[0.2cm]
Is it better for us to \textbf{Switch}, \textbf{Not Switch}, or is it \textbf{Equivalent}?
\\[0.2cm]
\textbf{Answer:} Not Switch.

\newpage
\section{Solutions}

\subsection{\hyperref[Es1]{Coins Under Cups}}
Let us define two sequences of random variables $\{O_i\}_{i=1}^{50},\{E_i\}_{i=1}^{50}$, representing respectively the result of the $i$-th toss of an odd-positioned coin and of an even-positioned coin, arranged on the table in sequential order, independent and identically distributed. It is clear that $\forall\;i:E_i\,,O_i \sim Ber(\frac{1}{2}).$ By assumption, Alice observes the outcomes of the coins in the order $(O_1, E_1, O_2, E_2, \ldots)$, while Bob observes the outcomes of the coins in the order $(O_1, O_2, \ldots, E_1, E_2, \ldots)$. The exercise asks to calculate:
\[
\displaystyle\frac{\mathbb P(\text{Alice wins})}{\mathbb P(\text{Bob wins})}
\]
We proceed as follows: fix a time $t$, calculate the probability that Alice wins exactly at time $t$, and then sum for $t=1,\dots,100$ to obtain $\mathbb P(\text{Alice wins})$. To do this, lets define this four events:
\begin{itemize}
    \item $T_1 = {\text{The game ends at time }t = 2m-1 ,\text{ for some } 1\leq m \leq 25};$
    \item $T_2 = {\text{The game ends at time }t = 2m-1 ,\text{ for some } 26\leq m \leq 50};$
    \item $T_3 = {\text{The game ends at time }t = 2m ,\text{ for some } 1\leq m \leq 25};$
    \item $T_4 = {\text{The game ends at time }t = 2m ,\text{ for some } 26\leq m \leq 49};$
\end{itemize}

Now, we calculate the probability that Alice wins in each of the previous events:

\begin{itemize}
    \item If $T_1$, Alice sees $O_m$ and wins if:
        \begin{enumerate}
            \item Bob has seen all tails in $(O_1,\dots,O_{2m-1})$ except for $O_m = H$ ;
            \item Alice has seen exactly one head among $(E_1,\dots,E_{m-1})$ ;
        \end{enumerate}
    that is, with probability:
    \[
    \mathbb P(\text{Alice wins} \cap T1) = \displaystyle\frac{m-1}{2^{3m-2}}.
    \]
    \item If $T_2$, Alice sees $O_m$ and wins if:
        \begin{enumerate}
            \item Bob has seen all tails in $(O_1,\dots,O_{50}, E_1,\dots,E_{2m-1-50})$ except for $O_m = H$ ;
            \item Alice has seen exactly one head among $(E_{2m-50},\dots,E_{m-1})$ ;
        \end{enumerate}
    that is, with probability:
    \[
    \mathbb P(\text{Alice wins}\cap T_2) = \displaystyle\frac{50-m}{2^{m+49}}.
    \]
    \item If $T_3$, the observed coins are $(O_1,\dots,O_{2m})$ and $(E_1,\dots,E_{m})$, for a total of $2^{3m}$ possible outcomes. Among these, we count in how many Alice wins.
        \begin{enumerate}
            \item Alice wins if $(O_1,\dots,O_{m})$ are all tails, if in $(O_{m+1},\dots,O_{2m})$ there is at most one head, and in $(E_1,\dots,E_{m-1})$ there is exactly one head, for a total of $(m+1)(m-1)$ favourable outcomes;
            \item Alice wins if there is exactly one head among $(O_1,\dots,O_{m})$, no head in $(O_{m+1},\dots,O_{2m})$, and no head in $(E_1,\dots,E_{m-1})$, for a total of $m$ favourable outcomes;
        \end{enumerate}
    that is, in this case Alice wins with probability:
    \[
    \mathbb P(\text{Alice wins}\cap T_3) = \displaystyle\frac{(m+1)(m-1)+m}{2^{3m}}= \displaystyle\frac{m^2+m-1}{2^{3m}}
    \]
    \item If $T_4$, the observed coins are $(O_1,\dots,O_{50})$ and $(E_1,\dots,E_{m})$, for a total of $2^{50+m}$ possible outcomes. Among these, we count in how many Alice wins. Define the following sets:
    \begin{itemize}
        \item  $C = (O_1,\dots,O_m,E_1,\dots,E_{2m-50})$, the coins seen by both Alice and Bob before $t$, $|C|= 3m-50$
        \item $Q = (E_{2m-49},\dots,E_{m-1})$, the coins seen only by Alice before $t$, $|Q|=49-m$;
        \item $R=(O_{m+1},\dots,O_{50})$, the coins seen only by Bob before $t$, $|R|=50-m$.
    \end{itemize}
    Then, in this case, Alice wins if:
        \begin{enumerate}
            \item There is no head in $C$, exactly one head in $Q$, and at most one head in $R$, for a total of $(49-m)(51-m)$ favourable cases;
            \item There is exactly one head in $C$, no head in $Q$, and no head in $R$, for a total of $3m-50$ favourable cases;
        \end{enumerate}
     that is, in this case Alice wins with probability:
    \[
    \mathbb P(\text{Alice wins}\cap T_4) = \displaystyle\frac{(49-m)(51-m)+3m-50}{2^{50+m}} = \displaystyle\frac{m^2 -97m +2449}{2^{50+m}}
    \]
    We observe that we exclude the case $m=50$, since neither player can win on the hundredth observation, given that both have seen the same coins.     
\end{itemize}
We can therefore state that:
\[
\mathbb P(\text{Alice wins})=
\displaystyle\sum_{m=1}^{25}\mathbb P(\text{Alice wins}\cap T_1) + 
\displaystyle\sum_{m=26}^{50}\mathbb P(\text{Alice wins}\cap T_2) + 
\]
\[
+\displaystyle\sum_{m=1}^{25}\mathbb P(\text{Alice wins}\cap T_3) +
\displaystyle\sum_{m=26}^{49}\mathbb P(\text{Alice wins}\cap T_4) \approx 0.3120
\]\\[0.2cm]
Now let us calculate, at the same fixed times and in the same way as before, Bob's probabilities of winning:
\begin{itemize}
    \item $ P(\text{Bob wins}\cap T_1)= \displaystyle\frac{m^2-m}{2^{3m-2}}$
    \item $P(\text{Bob wins}\cap T_2)= \displaystyle\frac{50-m}{2^{49+m}}$
    \item $P(\text{Bob wins}\cap T_3)= \displaystyle\frac{m^2+m-1}{2^{3m}}$
    \item $P(\text{Bob wins}\cap T_4)= \displaystyle\frac{50-m}{2^{50+m}}$
\end{itemize}
We can therefore state that:
\[
\mathbb P(\text{Bob wins})=
\displaystyle\sum_{m=1}^{25}\mathbb P(\text{Bob wins}\cap T_1) + 
\displaystyle\sum_{m=26}^{50}\mathbb P(\text{Bob wins}\cap T_2) + 
\]
\[
+\displaystyle\sum_{m=1}^{25}\mathbb P(\text{Bob wins}\cap T_3) +
\displaystyle\sum_{m=26}^{49}\mathbb P(\text{Bob wins}\cap T_4) \approx 0.4169 
\]\\[0.2cm]
and finally conclude that:
\[
\displaystyle\frac{\mathbb P(\text{Alice wins})}{\mathbb P(\text{Bob wins})} \approx \frac{0.3120}{0.4169 } \approx 0,75
\]
\\\(\qed\)

\subsection{\hyperref[Es2]{Head-Head or Head-Tails?}}
The exercise has already been analysed in a different and more general form by G. R. Grimmett in the article \textit{Alice and Bob on X: Reversal, Coupling, Renewal} \parencite{grimmett2025alicebobbbbx}. From the article one obtains the following two approximations for a generic number $n$ of coin tosses:
\[\begin{cases}
\mathbb P(\text{Bob wins}) - \mathbb P(\text{Alice wins}) = \displaystyle\frac{1}{2\sqrt{n\pi}}\\
\mathbb P(\text{Draw}) = \displaystyle\frac{1}{\sqrt{n\pi}}\implies \mathbb P(\text{Bob wins}) + \mathbb P(\text{Alice wins}) = 1-\frac{1}{\sqrt{n\pi}}
\end{cases}\]
Substituting $n=100$ and solving the linear system, we obtain that:
\[
\begin{cases}
    \mathbb P(\text{Alice wins}) \approx 0.458 \\
    \mathbb P(\text{Bob wins}) \approx 0.486
\end{cases}\implies\;\; \displaystyle\frac{\mathbb P(\text{Alice wins})}{\mathbb P(\text{Bob wins})} \approx 0.94
\]
\\\(\qed\)

\subsection{\hyperref[Es3]{A Sequence Race}}
We present the following state diagram, which represents the possible evolutions of the coin tosses in the sequence:

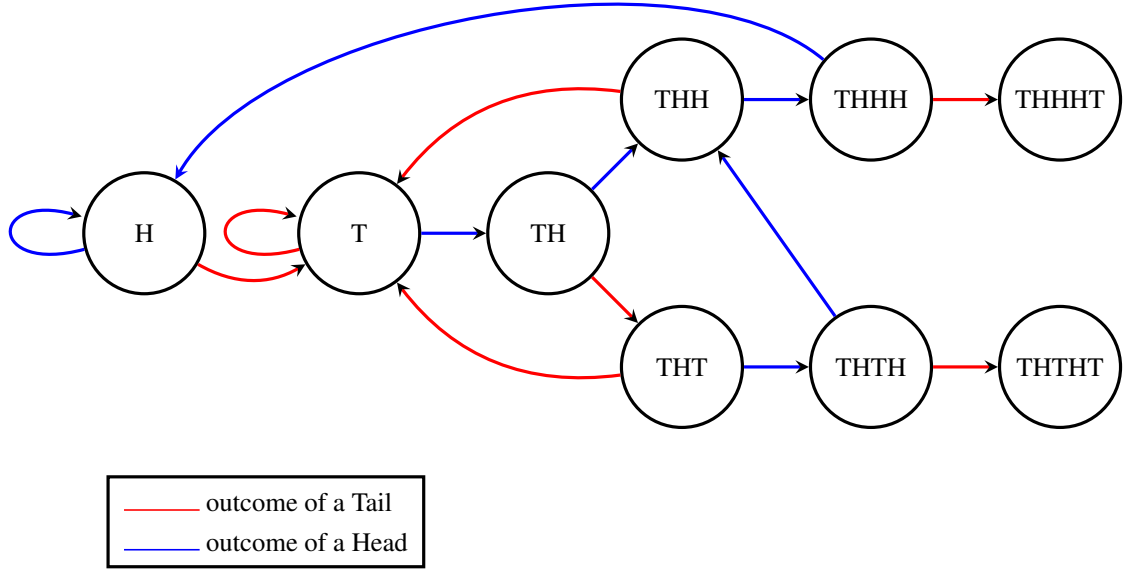
\begin{center}
\begin{tikzpicture}[
    ->,
    >=stealth,
    auto,
    node distance=2.5cm,
    very thick,
    state/.style={
        circle,
        draw,
        minimum size=1.6cm,
        inner sep=0pt,
        align=center
    }
]

    \node[state] (H) {H};
    \node[state, right=1.2cm of H] (T) {T};
    \node[state, right of=T] (TH) {TH};
    \node[state, below right of=TH] (THT) {THT};
    \node[state, above right of=TH] (THH) {THH};
    \node[state, right of=THT] (THTH) {THTH};
    \node[state, right of=THH] (THHH) {THHH};
    \node[state, right of=THTH] (THTHT) {THTHT};
    \node[state, right of=THHH] (THHHT) {THHHT};
    
    \draw[draw=red]
        (H) edge[bend right, below] node{} (T)
        (TH) edge[below left] node{} (THT)
        (THTH) edge[below] node{} (THTHT)
        (T) edge[loop left] node{} (T)
        (THT) edge[bend left, below] node{} (T)
        (THH) edge[bend right, above] node{} (T)
        (THHH) edge[above] node{} (THHHT);
    
    \draw[draw=blue]
        (TH) edge[above left] node{} (THH)
        (H) edge[loop left] node{} (H)
        (THT) edge[below] node{} (THTH)
        (THH) edge[above] node{} (THHH)
        (T) edge[above] node{} (TH)
        (THHH) to[out=140, in=60, looseness=0.7] node[above] {} (H)
        (THTH) edge[above right] node{} (THH);

    \node[draw, rectangle, align=left, anchor=north west, inner sep=6pt] 
        at (-0.5,-3.2) {
        \textcolor{red}{\rule{1cm}{0.6pt}} outcome of a Tail\\[4pt]
        \textcolor{blue}{\rule{1cm}{0.6pt}} outcome of a Head
        };
\end{tikzpicture}
\end{center}

By assumption, at every instant the probability of obtaining Head or Tail is the same, and the tosses are independent of one another. In this diagram the absorbing states are $THHHT$ and $THTHT$, which, with respect to the game, represent respectively Alice's victory and Bob's victory. Let $S$ be a generic state of the diagram; we then define:
\[p_S:= \mathbb P(\text{Alice wins}|\text{The current state is }S)\]
We observe that, by the law of total probability, we have:
\[p_H = \frac 12 p_H + \frac 12 p_T \implies p_H = p_S\]
Therefore, we shall set $\mathbb P(\text{Alice wins}) = p_H$. However, we also have that:
\[p_{T}=\frac 12 p_H + \frac 12 p_{TH} \implies p_{TH}= p_T = p_H\]
Similarly, we define all the other relations among the various states of the system, in order to form a system of equations:
\[
\begin{cases}
p_H=\frac12p_H+\frac12p_T,\\
p_T=\frac12p_T+\frac12p_{TH},\\
p_{TH}=\frac12p_{THH}+\frac12p_{THT},\\
p_{THH}=\frac12p_T+\frac12p_{THHH},\\
p_{THHH}=\frac12+\frac12p_H,\\
p_{THT}=\frac12p_T+\frac12p_{THTH},\\
p_{THTH}=\frac12p_{THH}.
\end{cases}
\]
Solving the system, we obtain that: 
\[\mathbb P(\text{Alice wins}) = p_H = \displaystyle\frac 59 \implies \mathbb P(\text{Bob wins}) = \frac 49\]
and therefore we conclude that:
\[\displaystyle\frac{\mathbb P(\text{Alice wins})}{\mathbb P(\text{Bob wins})} = \frac{5}{4} = 1.25\]
\\\(\qed\)

\subsection{\hyperref[Es4]{An Unfair Coin in a Haystack}}
Let X be the event ``the presenter has drawn the biased coin" and $H(i)$ the event ``the tossing of $i$ coins ten times each has produced only heads". Then, to answer the problem, we must calculate the conditional probability: $\mathbb{P}(X | H(9))$. Using Bayes' theorem, we obtain:
\begin{center}
    $\mathbb{P}(X | H(9)) = \displaystyle\frac{\mathbb{P}(H(9)|X) \mathbb{P}(X)}{\mathbb{P}(H(9))}$\\[0.2 cm]
\end{center}
But given that:
\begin{enumerate}
    \item $\mathbb{P}(H(9)|X) = (\frac{1}{{2}^{10}})^8$ 
    \item $\mathbb{P}(X)= 1 - (\frac{9999}{10000}\cdot\frac{9998}{9999}\cdot\frac{9997}{9998}\cdot\frac{9996}{9997}\cdot\frac{9995}{9996}\cdot\frac{9994}{9995}\cdot\frac{9993}{9994}\cdot\frac{9992}{9993}\cdot\frac{9991}{9992})= \frac{9}{10000}$
    \item $\mathbb{P}(H(9))= \mathbb{P}(H(9)|X)\cdot\mathbb{P}(X)+{P}(H(9)|X^c)\cdot\mathbb{P}(X^c)$
\end{enumerate}
We have:
\begin{center}
    $\mathbb{P}(X | T(9)) = \displaystyle\frac{(\frac{1}{1024})^8\frac{9}{10000}}{(\frac{1}{{2}^{10}})^8 \cdot \frac{9}{10000}+(\frac{1}{{2}^{10}})^9\cdot \frac{9991}{1000}} = 0.4798$
\end{center}
Thus Alice is right with a probability equal to $\sim 0.48$. 
\\ \( \qed \)

\subsection{\hyperref[Es5]{The Pen of Cows and Sheep}}
Let ${N}$ be the random variable describing the number of cows inside the new enclosure, ${N} \sim U([0,150])$. Let $FC$ be the event "the first animal to leave the new enclosure is a Cow" and $BC$ the event "both animals to leave the new enclosure are Cows". To answer the problem, we must therefore calculate: 
\begin{center}
    $\mathbb {P}(BC | FC) = \displaystyle\frac{\mathbb {P}(BC)}{\mathbb {P}(FC)}$\;.
\end{center}
We calculate the two values separately. The probability that the first animal is a Cow is:
\begin{center}
    $\mathbb {P}(FC)= \displaystyle\sum_{i=0}^{150} \mathbb {P}(FC|N=i) \mathbb {P}(N=i) = \displaystyle\frac{1}{151} \sum_{i=0}^{150} \frac{i}{150}=$ \\ $=\displaystyle\frac{1}{151 \cdot 150}\cdot \frac{151 \cdot 150}{2}=\frac{1}{2}$.
\end{center}
Again, the law of total probability is used to proceed the calculation of $\mathbb P(BC)$:
\begin{center}
$\mathbb{P}(BC) = \displaystyle\sum_{i=0}^{150} \mathbb{P}(BC|N=i) \mathbb{P}(N=i) = \displaystyle\frac{1}{151} \sum_{i=0}^{150} \frac{i (i-1)}{150\cdot 149} =$ \\[0.3 cm] $= \displaystyle\frac{1}{151\cdot150\cdot149}\cdot \left(\frac{150\cdot151\cdot301}{6}- \frac{150\cdot151}{2}\right)= \frac{1}{3}$\;.
\end{center}
To obtain this result we used the formula:
\begin{center}
    $\displaystyle\sum_{i=0}^{n}i^2 = \frac{n(n+1)(2n+1)}{6}$
\end{center}
Thus $\mathbb{P}(BC | FC) = \frac{2}{3}$\;, that is, knowing that the first animal to leave is a Cow, it is more likely that the second animal to leave is also a Cow. \\\( \qed \)

\newpage
\subsection{\hyperref[Es6]{Theft at the Greenhouse}}
Given the events:
\begin{itemize}
    \item 1P: "Alice saw that one of the two plants was a petunia";
    \item PO: "The man took a petunia and an orchid";
    \item PP: "The man took two petunias."
\end{itemize}
we calculate the probability that: $\mathbb P(PP|1P)$. Applying Bayes' theorem and the law of total probability, we obtain that:
\begin{center}
    $\mathbb P(PP|1P) = \displaystyle \frac{\mathbb P(PP)}{\mathbb P(1P)}=\frac{\displaystyle\frac{51}{101}\cdot\frac{50}{100}}{\mathbb P(1P|PP)\cdot \mathbb P(PP)+\mathbb P(1P|PO)\cdot \mathbb P(PO)}=$ \\[0.4cm] $\displaystyle\frac{\frac{51}{101}\cdot\frac{50}{100}}{\frac{51}{101}\cdot\frac{50}{100}+\frac{1}{2}\cdot(\frac{50\cdot 51}{101\cdot 100}+\frac{51\cdot 50}{100\cdot 101})} = \displaystyle\frac{1}{2}$.
\end{center}
The error one might fall into is confusing the event: "Bob knows that one specific plant is a petunia" with "Bob knows that at least one plant is a petunia".  \\\( \qed \)

\subsection{\hyperref[Es7]{Probabilistic Escape}}
\noindent Let $X = (X_1, \dots, X_{10})$ be the random vector representing the result of tossing ten fair independent coins, namely $X_i \sim \text{Ber}(\frac{1}{2})\;\;\; \forall i \in \{1, \dots, 10\}$. Then:
\begin{center}
    $\mathbb{P}(X = (x_1, \dots, x_{10})) = \left(\displaystyle\frac{1}{2}\right)^{10} \quad \forall x_i \in \{H, T\}\,, \quad i \in \{1, \dots, 10\}$.
\end{center}
From this it follows that the results reported by the two guards have the same probability of having actually occurred, and therefore both have the same chance of being the liar. However, assuming that:
\begin{itemize}
    \item the liar does not want to be discovered;
    \item the liar knows that Hercules might mistakenly think that the result reported by the first guard is less likely than the result reported by the second;
\end{itemize} then we could state that the liar is the second guard. Indeed, if the first were lying, why would he have chosen to report his result, which is apparently less credible in Hercules' logic, thereby leading Hercules to choose him as a potential liar? No: under the assumptions above, if the first were the liar, he would have chosen to report a result different from the one stated, and therefore it is more likely that he is telling the truth. Thus both "Equivalent" and "The second guard" are accepted as correct solutions.\\\( \qed \)

\subsection{\hyperref[Es8]{The Bounded Envelope Game Show}}
\noindent Define the following random variables: let $O$ be the value observed by Alice when opening the envelope chosen first, $M$ the maximum prize available for that evening, and $X$ the prize drawn. Finally, let $G$ be the final gain that Alice would obtain by choosing to switch envelopes; we want to calculate $\mathbb{E}[G]$. From here onwards we shall calculate everything in units for simplicity of notation, with the understanding that we are speaking of millions. By assumption, we have that $M \sim \text{Unif}\{10,\dots,25\}$, and that $X|M \sim \text{Unif}\{0,\dots,M\}$, whence, for $x \in \{0,\dots,25\}$:
\[
\mathbb{P}(X=x)=\displaystyle\sum_{m=10}^{25}\mathbb{P}(X=x|M=m)\mathbb{P}(M=m)= \frac{1}{16}\sum_{m=10}^{25}\mathbb{P}(X=x|M=m) = \frac{1}{16}\sum_{m=10}^{25}\displaystyle\frac{\mathds 1_{\{x\leq m\}} }{m}
\]
\noindent We can say that:
\[
\mathbb{P}(O=8)=\mathbb{P}(O=8|X=8)\mathbb{P}(X=8)+\mathbb{P}(O=8|X=16)\mathbb{P}(X=16)=\]
\[=\displaystyle\frac{1}{2}(\mathbb{P}(X=8)+\mathbb{P}(X=16)).
\]
\noindent Therefore, using Bayes' theorem, we can conclude that:
\[
\mathbb{E}[G]= 4\; \mathbb{P}(X=8|O=8) + 16\; \mathbb{P}(X=16|O=8) = 
\]
\[
= 4\cdot \displaystyle\frac{\mathbb{P}(O=8|X=8)\; \mathbb{P}(X=8)}{\mathbb{P}(O=8)}+16\cdot\frac{\mathbb{P}(O=8|X=16)\;\mathbb{P}(X=16)}{\mathbb{P}(O=8)}=
\]
\[
= \displaystyle\frac{4\;\mathbb{P}(X=8)+16\;\mathbb{P}(X=16)}{\mathbb{P}(X=8)+\mathbb{P}(X=16)}\approx\frac{0.7048}{0.0874}=8.0641.
\]
\\\( \qed \)

\subsection{\hyperref[Es9]{The Fruit Market}}
\noindent Let $R$ (resp. $G$) be the events "The blindfolded man took a red apple (a green apple)". Let $BG$ be the event "On the basis of his touch, the man believes he has taken a green apple". By assumption we have that: $\mathbb{P}(BG|G) = 0.85.$ Hence, using Bayes' theorem, we obtain that:
\[
\mathbb{P}(G|BG)=\displaystyle\frac{\mathbb{P}(BG|G)\cdot\mathbb{P}(V)}{\mathbb{P}(BG)}=\frac{\mathbb{P}(BG|G)\cdot\mathbb{P}(G)}{\mathbb{P}(BG|G)\cdot\mathbb{P}(G)+\mathbb{P}(BG|R)\cdot\mathbb{P}(R)}= 
\]
\[
= \displaystyle\frac{0.85 \cdot 0.1}{0.85\cdot 0.1 + 0,15 \cdot 0.9}=\frac{17}{44}.
\]
\\\( \qed \)

\subsection{\hyperref[Es10]{Holiday at the Sea}}
Consider a single family that decides to apply this strategy, and calculate the expected value of the difference between the number of boys and the number of girls among the children they will conceive. Represent the birth of a boy by \(+1\) and the birth of a girl
by \(-1\). Define the random variables:
\[
X_i =
\begin{cases}
+1, & \text{if the } i\text{-th child is a boy},\\
-1, & \text{if the } i\text{-th child is a girl}.
\end{cases}
\]
By assumption, the variables \(X_1,\dots,X_{10}\) are independent and identically
distributed, as Bernoulli random variables with parameter 1/2. Define:
\[
S_n=X_1+\cdots+X_n.
\]
Then \(S_n\) represents the difference between the number of boys and the number of
girls after \(n\) children. The family's strategy consists in stopping at the first instant at which this difference becomes positive, or in any case at the tenth child. We therefore define
the stopping time
\[
\tau=\min\{n\in\{1,\dots,10\}: S_n>0\},
\]
if such an \(n\) exists, and set $\tau=10$ in the case in which \(S_n\le 0\) for every \(n=1,\dots,10\). We observe that
\[
S_\tau=\sum_{k=1}^{10} X_k \mathbf 1_{\{\tau\ge k\}}.
\]
Indeed, the term \(X_k\) contributes to the final sum if and only if the family
gets to have at least \(k\) children, that is, if and only if \(\tau\ge k\). For every \(k\), the event \(\{\tau\ge k\}\) depends only on the first \(k-1\) children,
that is, on \(X_1,\dots,X_{k-1}\). Instead, \(X_k\) is independent of
\(X_1,\dots,X_{k-1}\). Therefore \(X_k\) is independent of
\(\mathbf 1_{\{\tau\ge k\}}\), and hence:
\[
\mathbb E\left[X_k\mathbf 1_{\{\tau\ge k\}}\right]
=
\mathbb E[X_k]\mathbb E\left[\mathbf 1_{\{\tau\ge k\}}\right]
\]
But since \(\mathbb E[X_k]=0\), we obtain:
\[
\mathbb E\left[X_k\mathbf 1_{\{\tau\ge k\}}\right]=0 \implies \mathbb E[S_\tau] = 0.
\]
Consequently, if the expected value of the difference between the number of male and female children is equal to zero for a generic family, it must be the case that the ratio of males to females in the population remains unchanged despite the application of the strategy, and is therefore 1:1.
\\\( \qed \)

\subsection{\hyperref[Es11]{Some Dice Games}}
The correct answer is that Carlo could have a greater chance of winning more games than everyone else, since there exists a configuration of dice assignable to each player that would satisfy the situation described in the exercise. Let us see it explicitly. Suppose that the dice of the respective players are:
\begin{itemize}
    \item Carlo: $[6,6,0,0,0,0]$
    \item Bob: $[3,3,3,3,3,3]$
    \item Alice: $[4,4,4,4,0,0]$
    \item Eva: $[5,5,1,1,1,1]$
\end{itemize}
Then, it is easily verified that, in the pairwise contests, the dice defined are consistent with the results of the victories:
\begin{itemize}
    \item $\displaystyle\frac{\mathbb{P}(\text{Eva wins against Alice})}{\mathbb{P}(\text{Alice wins against Eva})}= \frac{5
    }{9} \cdot \frac{9}{4} = \frac{5}{4}$
    \item $\displaystyle\frac{\mathbb{P}(\text{Alice wins against Bob})}{\mathbb{P}(\text{Bob wins against Alice})}= \frac{2
    }{3} \cdot \frac{3}{1} = 2$
    \item $\displaystyle\frac{\mathbb{P}(\text{Bob wins against Eva})}{\mathbb{P}(\text{Eva wins against Bob})}= \frac{2
    }{3} \cdot \frac{3}{1} = 2$
    \item $\displaystyle\frac{\mathbb{P}(\text{Eva wins against Carlo})}{\mathbb{P}(\text{Carlo wins against Eva})}= \frac{2
    }{3} \cdot \frac{3}{1} = 2$
    \item $\displaystyle\frac{\mathbb{P}(\text{Alice wins against Carlo})}{\mathbb{P}(\text{Carlo wins against Alice})}= \frac{4
    }{9} \cdot \frac{3}{1} = \frac{4}{3}$
    \item $\displaystyle\frac{\mathbb{P}(\text{Bob wins against Carlo})}{\mathbb{P}(\text{Carlo wins against Bob})}= \frac{2
    }{3} \cdot \frac{3}{1} = 2$
\end{itemize}
However, in the case in which all the players compete together, the chances are reversed in favour of Carlo; indeed:
\begin{itemize}
    \item $\mathbb{P}(\text{Eva wins against everyone})=\frac{2}{9}$
    \item $\mathbb{P}(\text{Alice wins against everyone})=\frac{8}{27}$
    \item $\mathbb{P}(\text{Bob wins against everyone})=\frac{4}{27}$
    \item $\mathbb{P}(\text{Carlo wins against everyone})=\frac{1}{3}$
\end{itemize}
Therefore, in this case, the dice configuration is consistent with the possibility that Carlo has a greater chance of winning more games than the others.
\\\( \qed \)

\subsection{\hyperref[Es12]{A Stochastic Triathlon Race}}

\subsection{\hyperref[Es13]{The New Broker}}
Let $A,B,F$ be the random variables indicating respectively the end-of-year return of Alberto, Bernardo, and Filippo. First of all, we calculate the expected values of each:
\begin{itemize}
    \item $\mathbb{E}[A]= 3\%$
    \item $\mathbb{E}[F]=2.96\%$
    \item $\mathbb{E}[B]= 3.32\%$
\end{itemize}
Therefore, it is clear that, in order to maximise profits, it is convenient to turn to Bernardo. Let us calculate the probability that Bernardo wins the investment fund prize, assuming that the returns of each broker are independent of one another:
\[
\mathbb{P}(\{B>F\} \cap \{B>A\})= \mathbb{P}(\{B=6, A=3, F=1\} \cup \{B=6, A=3, F=5\} \cup 
\]
\[\cup \{B=4, A = 3, F=1\})= 0.22\cdot 0.51 + 0.22\cdot 0.49 + 0.22 \cdot 0.51 = 0.3322\]
\\\( \qed \)

\subsection{\hyperref[Es14]{A Huge Monopoly Board}}

\subsection{\hyperref[Es15]{New Coins on the Ground}}
In the till there are $12$ coins, of which one in four is new; hence there are 3 new coins, while the non-new coins are $9$. The two coins given as change are chosen at random from the 12 coins in the till. Let $K$ be the number of new coins among the two received. Then $K$ has a hypergeometric distribution:
\[
\mathbb{P}(K=k)
=
\frac{\binom{3}{k}\binom{9}{2-k}}{\binom{12}{2}},
\qquad k=0,1,2.
\]
We therefore explicitly calculate its distribution:
\[
\mathbb{P}(K=0)=\frac{\binom{3}{0}\binom{9}{2}}{\binom{12}{2}}
=\frac{36}{66}=\frac{6}{11},
\]
\[
\mathbb{P}(K=1)=\frac{\binom{3}{1}\binom{9}{1}}{\binom{12}{2}}
=\frac{27}{66}=\frac{9}{22},
\]
\[
\mathbb{P}(K=2)=\frac{\binom{3}{2}\binom{9}{0}}{\binom{12}{2}}
=\frac{3}{66}=\frac{1}{22}.
\]
Now let $A$ be the event observed by the second man: \{at least one of the two coins that fell is new and landed heads\}. Denote by $B$ the event \{both coins landed heads\}. Then the required probability is
\[
\mathbb{P}(B\mid A)=\frac{\mathbb{P}(B\cap A)}{\mathbb{P}(A)}.
\]
First calculate $\mathbb{P}(A)$, conditioning with respect to the number $K$ of
new coins received. We have that:
\[
\mathbb{P}(A\mid K=0)=0,\qquad
\mathbb{P}(A\mid K=1)=\frac12,\qquad
\mathbb{P}(A\mid K=2)=\frac34.
\]
Indeed, if there are no new coins, the event $A$ is impossible; if there is only one, it must land heads; if there are two, at least one of the two must land heads.
Therefore, using the law of total probability, we obtain:
\[
\mathbb{P}(A) = 0\cdot \frac{6}{11} + \frac12\cdot \frac{9}{22} + \frac34\cdot \frac{1}{22} = \frac{21}{88}.
\]
Similarly, for the event $B\cap A$, we have that:
\[
\mathbb{P}(B\cap A\mid K=0)=0,\qquad
\mathbb{P}(B\cap A\mid K=1)=\frac14,\qquad
\mathbb{P}(B\cap A\mid K=2)=\frac14.
\]
Indeed, for $K=1$ or $K=2$, if both coins land heads, then certainly at least one new coin has landed heads. Thus:
\[
\mathbb{P}(B\cap A) = 0\cdot \frac{6}{11} + \frac14\cdot \frac{9}{22} + \frac14\cdot \frac{1}{22} = \frac{10}{88}.
\]
Finally, we can calculate the probability of the required event:
\[
\mathbb{P}(B\mid A) = \frac{\mathbb{P}(B\cap A)}{\mathbb{P}(A)} = \frac{\frac{10}{88}}{\frac{21}{88}} =
\frac{10}{21}.
\]
\\\( \qed \)

\subsection{\hyperref[Es16]{Chaos for the Kennel}}
Number the dogs from $0$ to $46$, where dog $0$ is the untrained one, while dogs $1,\dots,46$ are the trained ones. Denote by $C_i$ the kennel of dog $i$. Dog $0$ is released first and randomly chooses one among the $47$ available kennels. 
\begin{enumerate}
    \item If it chooses $C_0$, then all the trained dogs will find their own kennel free and in particular dog $46$ will enter $C_{46}$.
    \item If, instead, dog $0$ chooses kennel $C_k$, with $1\leq k\leq 46$, then dogs $1,\dots,k-1$ will certainly enter their respective kennels, while dog $k$, finding its own kennel occupied, will randomly choose a free kennel.
\end{enumerate}
Throughout the whole process, as long as neither of the two kennels $C_0$ and $C_{46}$ is chosen, every time a dog that finds its own kennel occupied, it chooses randomly among the kennels still free, including kennel 46. At that point:
\[
\begin{cases}
\text{if } C_0 \text{ is chosen, then dog } 46 \text{ will enter } C_{46},\\
\text{if } C_{46} \text{ is chosen, then dog } 46 \text{ will not enter } C_{46}.
\end{cases}
\]
Thus the final outcome depends only on which of the two kennels $C_0$ and
$C_{46}$ is chosen first by a dog forced to choose at random. By symmetry, these two possibilities have the same probability. Indeed, in every random choice among the kennels still available, the kennels $C_0$ and $C_{46}$, if both are still free, have the same probability of being chosen. Therefore:
\[
\mathbb{P}(\text{dog }46\text{ enters its own kennel})
=
\mathbb{P}(C_0 \text{ is chosen before } C_{46})
=
\frac12.
\]
\\\(\qed\)

\subsection{\hyperref[Es17]{Buses on Strike}}
Let $A$ be the event: \{No bus has passed for at least an hour\} and let $T$ be the random variable that measures how many minutes pass before the bus arrives. Therefore, the problem asks to calculate the quantity: $\mathbb E [T|A]$. Define the events $S_i$: \{there is a strike on the $i$-th day from today\}, $\forall\;i \in \mathbb N$, where for $i=0$ it is understood that \textit{there is a strike today}. In this way we can write that:
\[
\mathbb E [T|A] = \mathbb E[T|A,S_0^c]\cdot \mathbb P(S_0^c|A) +\mathbb E[T|A,S_0]  \cdot \mathbb P(S_0|A)\]
We calculate the various components separately. In the case in which there is no strike today, then $T|S_0^c\sim Exp(\frac{1}{10})$, hence:
\[\mathbb E[T|A,S_0^c]= 10\]
We use $\mathbb P (A|S_0^c) = \mathbb P (T>60|S_0^c)= e^{-6}$ and Bayes' theorem to calculate:
\[
\mathbb P(S_0^c|A) = \displaystyle\frac{\mathbb P (A|S_0^c) \cdot \mathbb P(S_0^c)}{\mathbb P (A|S_0^c) \cdot \mathbb P(S_0^c) + \mathbb P (A|S_0) \cdot \mathbb P(S_0)} \approx 0.1970
\]
Now calculate the conditional expected value in the case in which there is a strike today. Since it is $10$ in the morning, if there is a strike today I must certainly wait another 840 minutes before the beginning of the following day. Then define the random variable:
\[
N:=\min\{n\geq 1:S_n^c\},
\]
which indicates the first day, starting from tomorrow, on which there is no strike. Since the events $S_i$ are independent, we have that $N$ is a geometric random variable with parameter $0.99$, that is
\[
\mathbb P(N=n)=(0.01)^{n-1}\cdot 0.99,\qquad n\geq 1.
\]
Therefore:
\[
\mathbb E[N]=\frac{1}{0.99}.
\]
If $N=1$, then there is no strike tomorrow; I must wait $840$ minutes until midnight and then, since bus arrivals follow an exponential law with parameter $\lambda=\frac1{10}$, I must wait, on average, another $10$ minutes.
If instead $N>1$, then there are $N-1$ consecutive days of strike after today, and each contributes 1440 minutes of waiting. Thus, conditionally on $A$ and on $S_0$, we can write:
\[
T|A,S_0 \sim 840 + 1440(N-1) + X,
\]
where $X\sim Exp(\frac1{10})$ represents the waiting time, on the first day without a strike, before the arrival of the first bus. Therefore:
\[
\mathbb E[T\mid A,S_0]
=
840+1440\mathbb E[N-1]+ \mathbb E[X] = 840+1440\cdot\frac{0.01}{0.99}+10\approx 865
\]
It remains now to combine the two cases $S_0$ and $S_0^c$, and we conclude that:
\[
\mathbb E[T\mid A] = \mathbb E[T\mid A,S_0^c]\cdot \mathbb P(S_0^c\mid A) + \mathbb E[T\mid A,S_0]\cdot \mathbb P(S_0\mid A) \approx 696.
\]
\\\(\qed\)

\subsection{\hyperref[Es18]{Everyone's Wishes}}
Let $M,C,A$ be respectively the conditions desired by Mario, Camilla, and Aida. Let $\{X_i\}$ be the sequence of random variables representing the sex of the couple's children numbered in order of birth:
\[
X_i \sim Ber\left(\frac{1}{2}\right),\quad \forall\; i\in \mathbb N. 
\]
Let $VT$ be the event: "Everyone's wishes have been satisfied." and by assumption we know that $\exists\; T \in \mathbb N : (X_1,\dots,X_T)$ satisfies $VT$. The exercise asks us:
\[\mathbb P(X_2=M\; | \; VT) = \displaystyle \frac{\mathbb P(VT|X_2=M)\cdot \mathbb P (X_2=M)}{\mathbb P(VT)}\]
However, by assumption, we are also told that "the couple continues having children until $VT$ occurs"; hence from this we can infer that: $\mathbb P(VT)=1$, but also that $\mathbb P(VT|X_2=M)=1$. Moreover, we know that $\mathbb P(X_2=M)=\frac{1}{2}$, from which it follows that:
\[
\mathbb P(X_2=M\; | \; VT) = \displaystyle\frac{1}{2}
\]
\\\(\qed\)

\subsection{\hyperref[Es19]{The Dogs of Alice}}
Let $A$ be the event: "Hercules is a bad dog" and let $B$ be the event: "Hercules is a bad dog, but well cared for". It is evident that:
\[
B \subset A \implies \mathbb P(B) < \mathbb P(A)
\]
\\\(\qed\)

\subsection{\hyperref[Es20]{Monty Hall, But with a Golden Sheep}}
Let $A$ be the random variable indicating the content of the door initially chosen. Since the initial choice is random, we have:
\[
\mathbb{P}(A=\text{P})=\frac{1}{3}, \qquad
\mathbb{P}(A=\text{Pd'O})=\frac{1}{3}, \qquad
\mathbb{P}(A=\text{M})=\frac{1}{3}.
\]
where $P$ means a normal sheep, $Pd'O$ a golden sheep, and $M$ a car. Let $B$ be the random variable indicating the content of the door opened by the presenter. We observe that in no scenario will the presenter open the door containing the car. Therefore, a priori:
\[
\mathbb{P}(B=\text{P})=\frac{1}{2}, \qquad
\mathbb{P}(B=\text{Pd'O})=\frac{1}{2}, \qquad
\mathbb{P}(B=\text{M})=0.
\]
We want to calculate the probability that the door initially chosen contains the golden sheep, knowing that the presenter has opened a door containing a sheep. Applying Bayes' theorem, we obtain:
\[
\mathbb{P}(A=\text{Pd'O} \mid B=\text{P}) = \frac{ \mathbb{P}(B=\text{P} \mid A=\text{Pd'O})\cdot \mathbb P (A=\text{Pd'O})}{\mathbb{P}(B=\text{P})}.
\]
If the door initially chosen contains the golden sheep, then the presenter, being unable to open either the door with the car or the one chosen by the contestant, is forced to open the door containing the sheep. Thus:
\[
\mathbb{P}(B=\text{P} \mid A=\text{Pd'O})=1.
\]
Therefore:
\[
\mathbb{P}(A=\text{golden sheep} \mid B=\text{sheep})=\frac{1 \cdot \frac{1}{3}}{\frac{1}{2}}=\frac{2}{3}.
\]
from which it follows that it is convenient \textbf{Not to Switch} one's initial choice.
\\\(\qed\)

\newpage
\printbibliography

\end{document}